\documentclass[12pt]{article}
\usepackage{epsfig}
\usepackage{algorithm,algorithmic}
\usepackage{amsmath}
\usepackage{amssymb}
\usepackage{url}

\usepackage{times}
\usepackage{graphicx}
\usepackage{latexsym}
\usepackage{amsfonts}
\usepackage{enumerate}
\usepackage{tabularx}

\begin{document}
\title{\vspace{-17mm}
FFT-based Alignment of 2d  Closed Curves
with Application to Elastic Shape Analysis}
\author{\tt\small G\"unay Dogan$^{1,2}$, Javier Bernal$^2$, Charles Hagwood$^2$\\
$^1${\small \sl Theiss Research,} \\
{\small \sl La Jolla, CA, USA} \\
$^2${\small \sl National Institute of Standards and Technology,} \\
{\small \sl Gaithersburg, MD 20899, USA} \\
{\tt\small $\{$gunay.dogan,javier.bernal,charles.hagwood$\}$} \\
{\tt\small @nist.gov}}
\date{\ }
\maketitle
\vspace{-13mm}
\begin{abstract}
For many shape analysis problems in computer vision and scientific imaging
(\emph{e.g.}, computational anatomy, morphological cytometry),
the ability to align two closed curves in the plane is crucial.
In this paper, we concentrate on rigidly aligning pairs of closed curves
in the plane. If the curves have the same length and are centered
at the origin, the critical steps to an optimal rigid alignment are finding the
best rotation for one curve to match the other and redefining the starting
point of the rotated curve so that the starting points of the two curves match.
Unlike open curves, closed curves do not have fixed starting points, and this
introduces an additional degree of freedom in the alignment.
Hence the common naive method to find the best rotation and starting point
for optimal rigid alignment has $O(N^2)$ time complexity, $N$ the
number of nodes per curve. This can be slow for curves with large numbers
of nodes. In this paper, we propose a new $O(N \log N)$ algorithm for this
problem based on the Fast Fourier Transform. Together with uniform
resampling of the curves with respect to arc length, the new algorithm
results in an order of magnitude speed-up in our experiments.
Additionally, we describe how we can use our new algorithm as part of
elastic shape distance computations between closed curves to obtain
accurate shape distance values at a fraction of the
cost of previous~approaches.
\end{abstract}

\setcounter{page}{1}

\section{Introduction}\label{S:intro}

For many problems in computer vision and scientific imaging,
one needs to quantitatively
compare boundaries of objects, \emph{e.g.}, organs in computational
anatomy, cells in morphological cytometry. Often, the comparison
of the boundaries of two objects requires that they be aligned first.
Thus, methods for aligning boundaries of 2d objects have been discussed in
\cite{ayache, cohen, larsen, li, schwartz, sebastian, umeyama}.
Mathematically, the problem is that of aligning pairs of 2d closed curves.
In this article, we concentrate on rigidly aligning pairs of 2d closed curves,
assumed to be square-integrable.
The curves have been processed to match in scale (say, by normalization
to unit length) and are centered at the origin of a reference coordinate
system. In practice, curves are given as lists of points,
each list usually in a counterclockwise order around its curve,
with the choice of the first point in a list usually arbitrary.
With curves in this form, optimal rigid alignment is achieved by\\
1) \emph{Rotating} one curve, \emph{i.e.}, points that define curve, to match
as much as possible the points that define the other curve.\\
2) Among points that define the rotated curve, choosing the best
\emph{starting point} to match the starting point of the other~curve.

In this paper, we focus on these two factors and propose a new algorithm
to find the optimal pair of rotation and starting point.
If we do not account for these two factors, we may fail to align
even two different versions of the \emph{same} curve.
However an additional difficulty could prevent our method
from working properly: the discrete sampling of one curve
and the placing of the representative nodes often will be incompatible
with those of the other, especially if they come from different measured
data.  As will be described below, this difficulty can usually be alleviated
by uniform resampling of the two curves with respect to arc length.

If the starting point is known, the optimal
rotation can be found using the Kabsch algorithm \cite{kabsch1,kabsch2} in
$O(N)$ time, $N$ the number of nodes per curve (see Section~\ref{S:optim}).
Otherwise, a common approach is to loop over all starting point candidates
while computing optimal rotation for each, and then choosing the pair of
starting point and rotation with optimal alignment. This approach is
$O(N^2)$ and can be slow for curves with large numbers
of nodes. In this paper, we propose a new fast $O(N\log N)$ method for
this problem based~on the Fast Fourier Transform (FFT), used
together with uniform resampling of the curves with respect to arc length.
We then use this method as part of elastic shape distance computations
between closed curves and find that we compute accurate shape distances in
significantly less time than previous approaches
\cite{dogan,srivastava2}.

\section{Optimal starting point and rotation}\label{S:optim}

In this section, we describe our new algorithm to compute optimal starting
point and \mbox{rotation} for rigid alignment of closed curves. Initially we
have two curves (not necessarily closed) of unit length, centered at the origin,
defined by functions $\beta_i: [0,1]\rightarrow \mathbb{R}^2$, $i=1,2$,
assumed to be square-integrable, \emph{i.e.},
$\|\beta_i\|_{L^2}^2 = \int_0^1 \|\beta_i(t)\|^2 dt < \infty$, $i=1,2$.
First we review how the optimal rotation
$
R=
R(\theta) = \left(
\begin{smallmatrix}
\cos(\theta) & \sin(\theta) \\
-\sin(\theta) & \cos(\theta)\\
\end{smallmatrix} \right)
$
is found through angle $\theta$ that minimizes mismatch~energy
\begin{equation}\label{E:energy0}
E_0(R) = \int_0^1 \| \beta_1(t) - R \beta_2(t) \|^2 dt.
\end{equation}
%
For this purpose we first rewrite energy~\eqref{E:energy0} as follows
\begin{equation*}
E_0(R) = \|\beta_1\|_{L^2}^2 + \|\beta_2\|_{L^2}^2
-2 \int_0^1 \beta_1^T(t) R \beta_2(t) dt,
\end{equation*}
where $\|\beta_i\|_{L^2}^2$, $i=1,2$, have constant value
(note $\|R \beta_2(t)\| = \|\beta_2(t)\|$, since $R$ is a rotation matrix).
Then minimizing energy~\eqref{E:energy0} is equivalent to maximizing
%
\begin{equation*}
\int_0^1 \beta_1^T(t) R \beta_2(t) dt = tr( R A^T ),
\end{equation*}
where $tr(R A^T )$ is the trace of $R A^T$ and $A$ is the $2\times 2$ matrix
defined by $ A_{kj} = \int_0^1 \beta_{1k}(t) \beta_{2j}(t) dt$, $k,j=1,2.$
By the Singular Value Decomposition (SVD) theorem $A=USV^T$, where $U$, $V$
are  $2\times 2$ orthogonal matrices,
$S = \left(\begin{smallmatrix} \sigma_1 & 0 \\
                                      0 & \sigma_2 \end{smallmatrix}\right)$,
$\sigma_1 \geq \sigma_2\geq 0$. Then $tr(RA^T)= tr(RVSU^T)= tr(SU^TRV) = tr(ST)$,
where $T=U^TRV$.
With $T = \left(\begin{smallmatrix} t_{11} & t_{12} \\
                                    t_{21} & t_{22} \end{smallmatrix}\right)$,
then $tr(ST)=\sigma_1t_{11} + \sigma_2t_{22}$. Since $U$, $R$, $V$ are orthogonal
so is $T$. Thus $-1\leq t_{kj} \leq 1$, $k,j=1,2$. If $\det(U)\det(V)>0$ then
$\det(T)=1$ so that the maximum value of $tr(ST)$ occurs when $t_{11}= t_{22}=1$ which
implies $T=I$, the identity matrix. Thus $U^TRV=T=I$ and $R=UV^T$ is an optimal
rotation. If $\det(U)\det(V)<0$ then $\det(T)=-1$ which implies
$t_{11}t_{22}-t_{21}t_{12} = -1$ so that $t_{11}t_{22}\leq 0$. Since $|t_{11}|=|t_{22}|$
then the maximum value of $tr(ST)$ occurs when $t_{11}=1$, $t_{22}=-1$ which implies
$T =\tilde{S} = \left(\begin{smallmatrix} 1 & 0 \\
                                           0 &-1 \end{smallmatrix}\right)$.
Thus $U^TRV=T=\tilde{S}$ and $R=U\tilde{S}V^T$ is an optimal rotation.

In the case of \emph{closed} curves, we have two such curves
of unit length, centered at the origin, defined by periodic functions
$\beta_i:\mathbb{R} \rightarrow \mathbb{R}^2, \ \beta_i(t+1)=\beta_i(t)$ for all
values of $t$, $i=1,2$, assumed to be square-integrable over any finite interval
of~$\mathbb{R}$. In practice $\beta_1$ and $\beta_2$ are given as finite lists of
points, say $N$ points per curve for some integer $N>0$, each list in a
counterclockwise order around its curve with first and last points the same.
At first we parametrize these curves with the discrete uniform parametrization of
interval $[0,1]$ obtained from partitioning $[0,1]$ into subintervals with
endpoints $t_l=(l-1)h$, $h=1/(N-1)$, $l=1,\ldots,N$, by defining
$\beta_i^l = (\beta_{i1}^l,\beta_{i2}^l)^T$ by $\beta_i^l = \beta_i(t_l)$, $i=1,2$,
$l=1,\ldots,N$, where $\beta_i(t_l)$ is the $l^{th}$ point in the list for~$\beta_i$.
We then take advantage of the curves being uniformly parametrized this way to discretize
integral~\eqref{E:energy0} using the uniform trapezoidal rule for closed~curves:
\begin{equation}\label{E:disc-energy0}
E^h_0(R) = h \sum_{l=1}^{N-1} \| \beta_1^l - R \beta_2^l \|^2.
\end{equation}
As we did with \eqref{E:energy0}, we can rewrite
\eqref{E:disc-energy0} as follows
\begin{equation*}
E^h_0(R) = h(\sum_{l=1}^{N-1} (\|\beta_1^l\|^2 + \|\beta_2^l\|^2)
- 2(\sum_{l=1}^{N-1} (\beta_1^l)^T R \beta_2^l)),
\end{equation*}
so that minimizing \eqref{E:disc-energy0} is equivalent to maximizing
\begin{equation*}
\sum_{l=1}^{N-1} (\beta_1^l)^T R \beta_2^l = tr(RA^T),
\end{equation*}
where $A$ is the $2\times 2$ matrix defined by
$A_{kj} = \sum_{l=1}^{N-1} \beta_{1k}^l\beta_{2j}^l$, $k,j = 1,2$.
Then as we did for~\eqref{E:energy0},
optimal rotation $R$ for~\eqref{E:disc-energy0} can be computed
with the SVD of $A$ or, more precisely, with the Kabsch algorithm
(see Algorithm~\ref{A:kabsch}) \cite{kabsch1,kabsch2}.

\begin{algorithm}
\caption{Computing optimal $R$ for starting point $t_0 = 0$}
\label{A:kabsch}
\begin{algorithmic}
\STATE Compute $A_{kj} = \sum_{l=1}^{N-1}\beta_{1k}^l \beta_{2j}^l, \ k,j=1,2.$
\STATE Compute SVD of $A$ s.t.\! $A = U S V^T$.
\STATE If $\det(U)\det(V) > 0$
       then $\tilde{S} = \left(\begin{smallmatrix} 1 & 0 \\
                                                   0 & 1 \end{smallmatrix}\right)$
       else $\tilde{S} = \left(\begin{smallmatrix} 1 & 0 \\
                                                   0 &-1 \end{smallmatrix}\right)$.
\STATE Return $R = U \tilde{S} V^T$.
\end{algorithmic}
\end{algorithm}

Now we find  starting point $t_0$ and rotation $R$ that
give optimal rigid alignment of $\beta_1$ and $\beta_2$ by minimizing
mismatch energy
\begin{equation}\label{E:energy}
E(t_0,R) = \int_0^1 \| \beta_1(t) - R \beta_2(t+t_0) \|^2 dt,
\quad t_0 \in [0,1].
\end{equation}

Again with $\beta_1$ and $\beta_2$ uniformly parametrized as above,
for each $m$, $1\leq m\leq N-1$, we define
$\beta_2^{l\oplus m} = (\beta_{21}^{l\oplus m},\beta_{22}^{l\oplus m})^T$, by
\begin{equation*}
\beta_2^{l\oplus m} = \beta_2(t_l+t_m),\ l=1,\ldots,N.
\end{equation*}
With $t_0= t_m$, we then discretize integral~\eqref{E:energy}
using again the uniform trapezoidal rule for closed curves:
\begin{equation}\label{E:disc-energy}
E^h(t_0,R) = h \sum_{l=1}^{N-1} \| \beta_1^l - R \beta_2^{l\oplus m} \|^2.
\end{equation}
In addition, for each $m$, we define $2\times 2$ matrix $A(t_m)$ by
%
%
\begin{equation}\label{E:A-mtx}
A_{kj}(t_m) = \sum_{l=1}^{N-1} \beta_{1k}^l\beta_{2j}^{l\oplus m}, k,j = 1,2.
\end{equation}
For $t_0=t_1,\ldots,t_{N-1}$,
we could use Algorithm~\ref{A:kabsch} to compute $A=A(t_0)$ and the best
rotation $R(t_0) = R$, and then return the pair $(t_0,R(t_0))$
that gives the highest value for $tr(RA^T)$, \emph{i.e.}, the lowest value
for~\eqref{E:disc-energy}. This commonly used approach is $O(N^2)$
as the computation of $A(t_0)$ for a single $t_0$ is $O(N)$.

With $A_{kj}(t_m)$ as in~\eqref{E:A-mtx}, we propose to compute vectors
$A_{kj}=$\\ $(A_{kj}(t_1),\ldots,A_{kj}(t_{N-1}))$, $k,j=1,2$,
in $O(N\log N)$ time using
FFT to accomplish the Discrete Fourier Transform (DFT).
For this purpose we define $\tilde{\beta}_1^l =
(\tilde{\beta}_{11}^l,\tilde{\beta}_{12}^l)^T$ by
$\tilde{\beta}_1^l = \beta_1^{N-l+1}$, $l=1,\ldots,N$, and vectors
$\tilde{\beta}_{1k}$ and $\beta_{2j}$ by
$\tilde{\beta}_{1k}=(\tilde{\beta}_{1k}^1,\ldots,\tilde{\beta}_{1k}^{N-1})$,
$\beta_{2j}=(\beta_{2j}^1,\ldots,\beta_{2j}^{N-1})$, $k,j=1,2$.
Given arbitrary vectors $x$, $y$ of length~$N-1$, we let {\bf DFT}$(x)$ and
{\bf DFT}$^{-1}(y)$ denote the DFT of~$x$ and the inverse DFT of~$y$, respectively.
With the symbol $\cdot$ indicating component by component multiplication of two vectors,
then by the convolution theorem for the DFT, we have for $k,j = 1,2,$
\begin{equation*}
A_{kj}=
(\sum_{l=1}^{N-1} \beta_{1k}^l\beta_{2j}^{l\oplus 1},\ldots,
 \sum_{l=1}^{N-1} \beta_{1k}^l\beta_{2j}^{l\oplus (N-1)})=
\mathrm{\bf DFT}^{-1}[\mathrm{\bf DFT}(\tilde{\beta}_{1k})\cdot
\mathrm{\bf DFT}(\beta_{2j})],
\end{equation*}
%
which enables us to reduce the computation of the matrix
element $A_{kj}(t_0)$ for all $t_0$ to three $O(N\log N)$ FFT operations
(thus a total of twelve for all of $A(t_0)$).
Once we have computed $A(t_0)$ for all $t_0$ this way, we can loop over
$t_0$ candidates, compute for each $t_0$ the corresponding optimal rotation $R(t_0)$
and then $tr(R(t_0)A(t_0)^T)$ (instead of $E^h(t_0,R(t_0))$),
and return the pair $(t_0,R(t_0))$ that gives the highest $tr(R(t_0)A(t_0)^T)$
value, \emph{i.e.}, the lowest value for~$E^h(t_0,R(t_0))$ in~\eqref{E:disc-energy}.
This is summarized in Algorithm~\ref{A:fast-optim}.
There for arbitrary vectors $x$, $y$ of length~$N-1$, {\bf FFT}$(x)$, {\bf IFFT}$(y)$
denote {\bf DFT}$(x)$, {\bf DFT}$^{-1}(y)$, respectively, computed with FFT.
Note that the computation of the SVD of $A(t_0)$ has $O(1)$ complexity.
In addition, the computation of $tr(R(t_0)A(t_0)^T)$ has $O(1)$ complexity
as opposed to that of $E^h(t_0,R(t_0))$ in~\eqref{E:disc-energy}
that would have $O(N)$ complexity.

\begin{algorithm}
\caption{Fast algorithm for optimal $(t_0, R(t_0))$}
\label{A:fast-optim}
\begin{algorithmic}
\STATE Compute $A_{kj}=(A_{kj}(t_1),\ldots,A_{kj}(t_{N-1}))=$
              {\bf IFFT}[{\bf FFT} $(\tilde{\beta}_{1k})\,\cdot\,$
              {\bf FFT}$(\beta_{2j})$], $k,j=1,2$.
\FOR {$t_0 = t_1, \ldots, t_{N-1}$}
\STATE Compute SVD of $A = A(t_0)$ s.t.\! $A = U S V^T$.
\STATE If $\det(U)\det(V) > 0$
       then $\tilde{S} = \left(\begin{smallmatrix} 1 & 0 \\
                                                   0 & 1 \end{smallmatrix}\right)$
       else $\tilde{S} = \left(\begin{smallmatrix} 1 & 0 \\
                                                   0 &-1 \end{smallmatrix}\right)$.
\STATE Set $R = R(t_0) = U \tilde{S} V^T$.
\STATE Compute $tr(RA^T)$.
\ENDFOR
\STATE Return $(t_0,R(t_0))$ that gives the highest $tr(RA^T)$.
\end{algorithmic}
\end{algorithm}

%
%
Above we described the obvious way of parametrizing the curves from their
lists of points with the same uniform parametrization of~$[0,1]$.
It should be pointed out that for the development of the two algorithms
above, parametrizing the curves uniformly in any other way would have
produced the same results. As noted in the introduction, the discrete
sampling of one curve and the placing of the representative nodes may not
be compatible with those of the other. Accordingly, when rigidly aligning curves
with the algorithms above, the obvious way of parametrizing the curves
uniformly may not be the most appropriate. Thus, for this purpose,
it seems more natural to parametrize each curve by its arc length with the
same uniform parametrization of~$[0,1]$. In particular, in the case of similar
curves, optimal $(t_0,R)$ values are then produced, as the parametrized curves
then match in a way enabling correspondence
between the curves at endpoints of subintervals in the parametrization
domain~$[0,1]$.
We do this in two steps. In the first step, for $i=1,2$, we define
discrete parameter~$s^l_i$~by 
%
\begin{equation}\label{E:arc-length-param}
s^l_i = L^l_i/L^N_i, l = 1,\ldots,N,
\quad
L^l_i = \sum_{m=2}^{l} \|\beta^m_i - \beta^{m-1}_i\|, \ l=2,\ldots,N,
\quad L^1_i = 0.
\end{equation}
%
In the second step, we strive to approximately parametrize $\beta_1$ and $\beta_2$
by their arc lengths with the same uniform parametrization of~$[0,1]$
using $s_1^l$ and $s_2^l$ in~\eqref{E:arc-length-param}.
As we already have parameter and curve node pairs
$\{(s_1^l,\beta_1^l)\}_{l=1}^N$ for $\beta_1$,
$\{(s_2^l,\beta_2^l)\}_{l=1}^N$ for $\beta_2$,
we can interpolate the curves with cubic splines using these parameter and curve
node pairs, and obtain new curve nodes $\beta_i^l$, $i=1,2$, $l=1,\ldots,N$,
by evaluating the cubic spline interpolants at $t_l=(l-1)h$, $h=1/(N-1)$,
$l=1,\ldots,N$, thus approximately parametrizing each curve by its
arc length with the same uniform parametrization of~$[0,1]$.
This resampling procedure may introduce some loss in the accuracy of
the curve representation (effectively a geometric approximation
error), which will have a minor impact on the optimal $(t_0,R)$
computed. We find that the error introduced by resampling is not
important for our main application, the computation of elastic shape
distances. The algorithm we use for shape distance computation
is the iterative algorithm in~\cite{dogan} which
improves on an imperfect estimate of $(t_0,R)$. We could of course
use nonuniform parametrizations of~$[0,1]$ to parametrize the curves
by their arc lengths, such as those based on curvature~\cite{cui}.
However in that case, a nonuniform FFT would be required.

\begin{figure}
\begin{center}
\begin{tabular}{|c|c|c|}
\hline
\includegraphics[width=.1\textwidth]{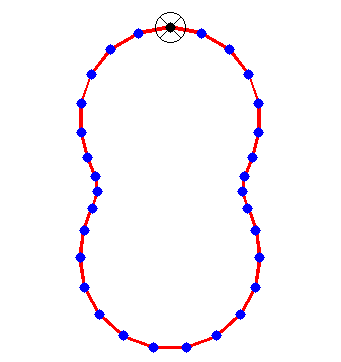}
&
\includegraphics[width=.1\textwidth]{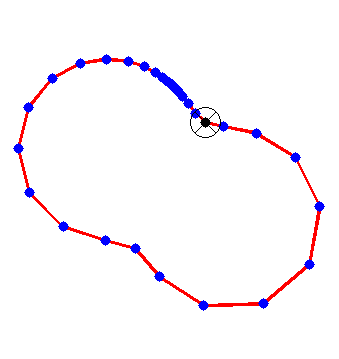}
&
\includegraphics[width=.1\textwidth]{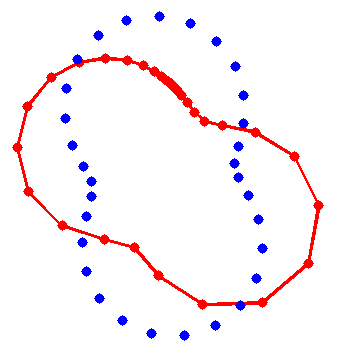}
\\ \hline
\includegraphics[width=.1\textwidth]{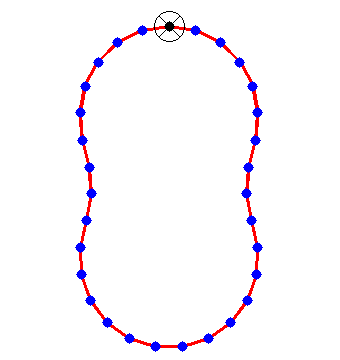}
&
\includegraphics[width=.1\textwidth]{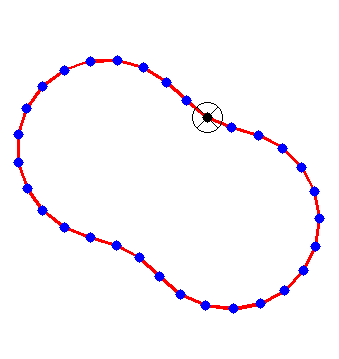}
&
\includegraphics[width=.1\textwidth]{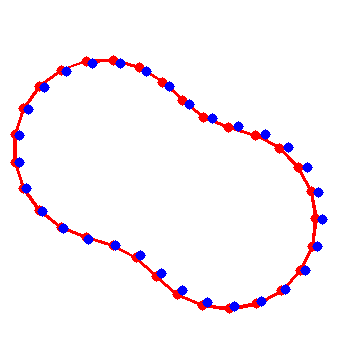}
\\ \hline
\end{tabular}
\end{center}
\caption{\label{F:resampling}
Impact of sampling and node placement in aligning
the curves in the 1st and 2nd columns.
When the sampling of nodes in the curve is non-uniform
(2nd in top row), the alignment fails (3rd in top row).
When the sampling is relatively uniform (bottom row), the points
in one curve can be easily matched with corresponding nodes
in the other curve to ensure successful alignment.
}
\end{figure}

In Figure~\ref{F:resampling}, using a hippopede curve, we illustrate
the impact of sampling and node placement on the computation of optimal
$(t_0,R)$. Top row of Figure~\ref{F:resampling}
includes two different samplings of the curve, one relatively uniform
(left-most column) and another one (middle column), which has nodes
concentrated in one part of the curve. We are not able to align the
two curves well (right-most column) by simply parametrizing them in
the obvious way that does not take arc length into consideration.
The figure includes resampled versions of
the two curves (in the bottom row, same order). As the resampled
curves have the same uniform arc length parametrization that
ensures good correspondence between the curve nodes, we are able to
align the two curves very well. In the figure the first points in
the lists of points (initial starting points) are circled in black.
The angle of rotation~is~$\frac{\pi}{3}$.

\begin{figure}
\begin{center}
\begin{tabular}{|ccccc|}
\hline
\includegraphics[width=.0375\textwidth]{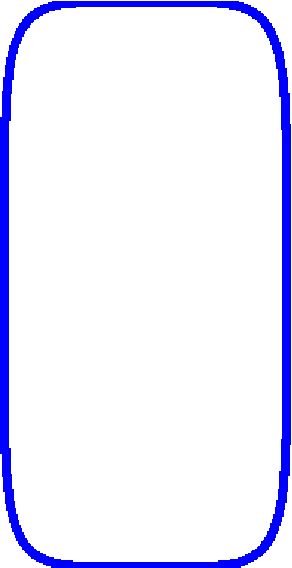}
&
\includegraphics[width=.05\textwidth]{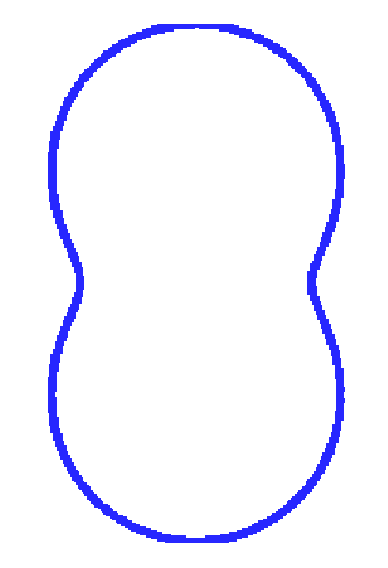}
&
\includegraphics[width=.0675\textwidth]{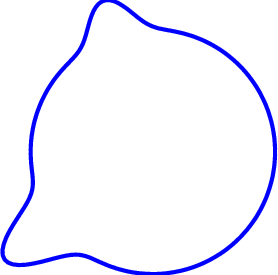}
&
\includegraphics[width=.0675\textwidth]{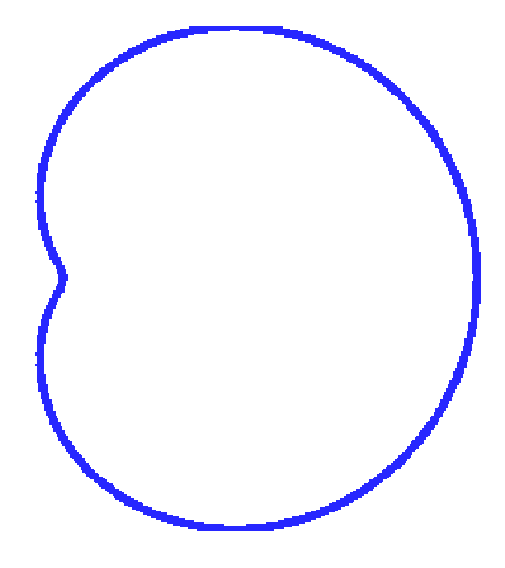}
&
\includegraphics[width=.075\textwidth]{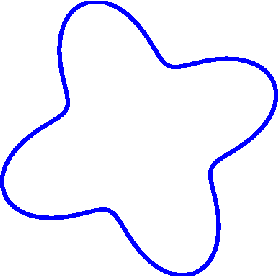}
\\ \hline
\includegraphics[height=.075\textwidth]{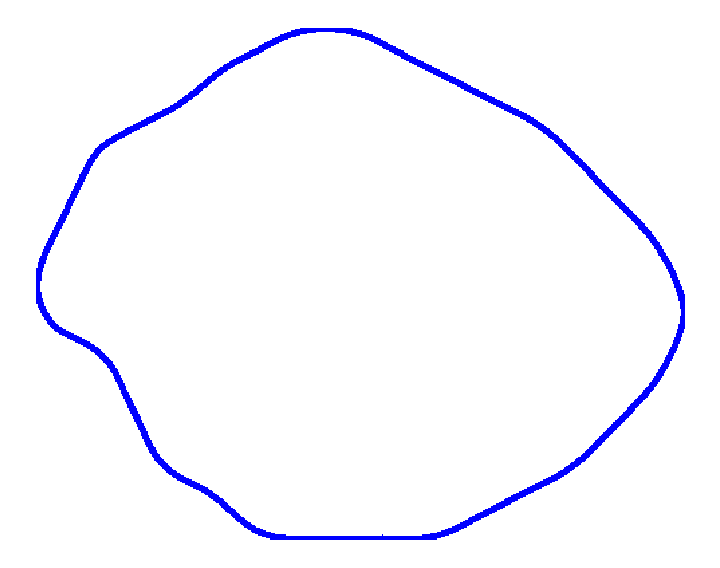}
&
\includegraphics[height=.075\textwidth]{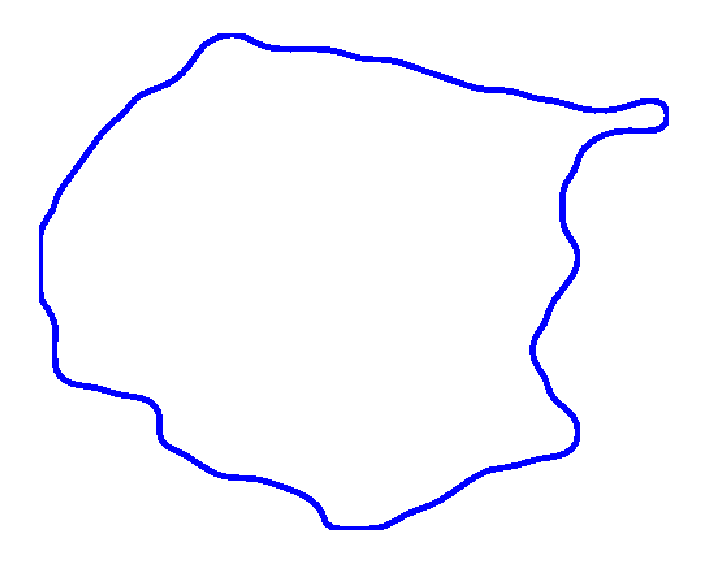}
&
\includegraphics[height=.075\textwidth]{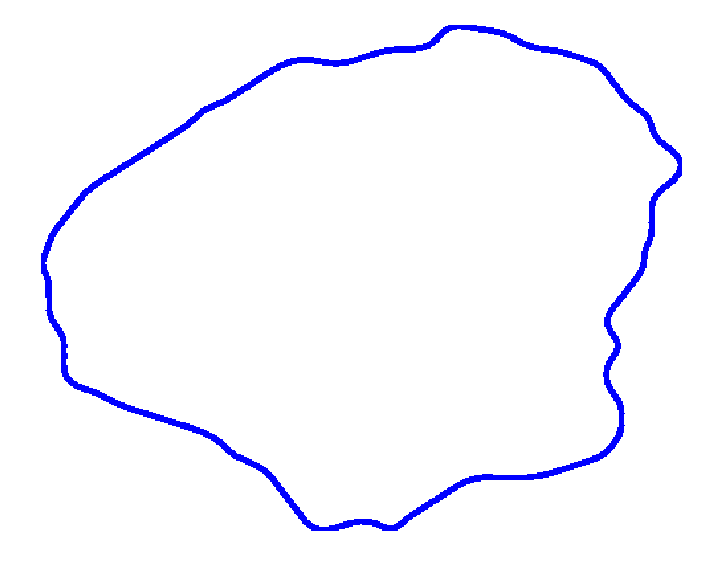}
&
\includegraphics[height=.075\textwidth]{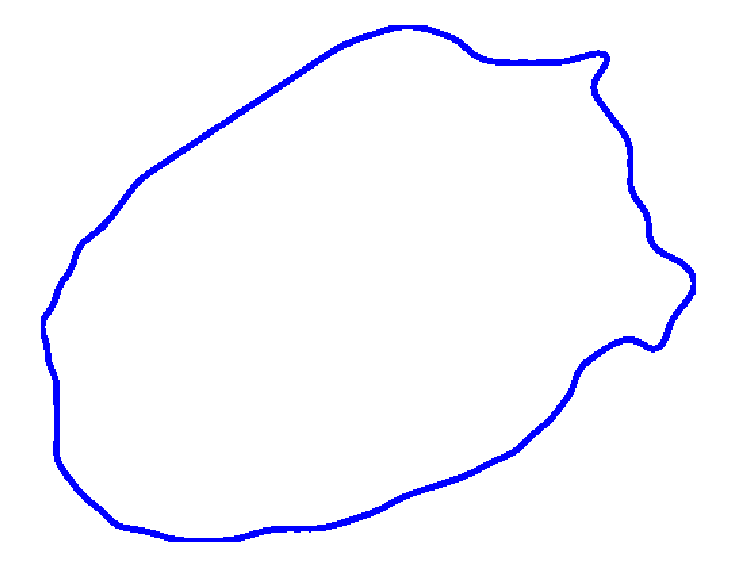}
&
\includegraphics[height=.075\textwidth]{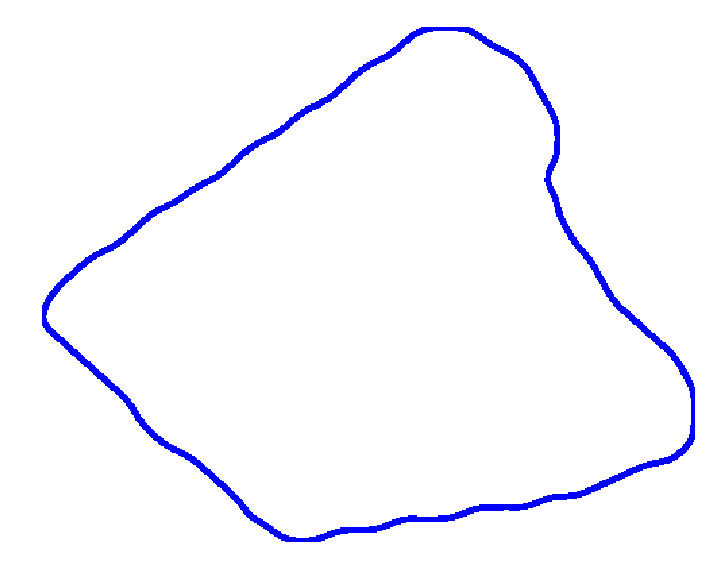}
\\ \hline
\includegraphics[height=.075\textwidth]{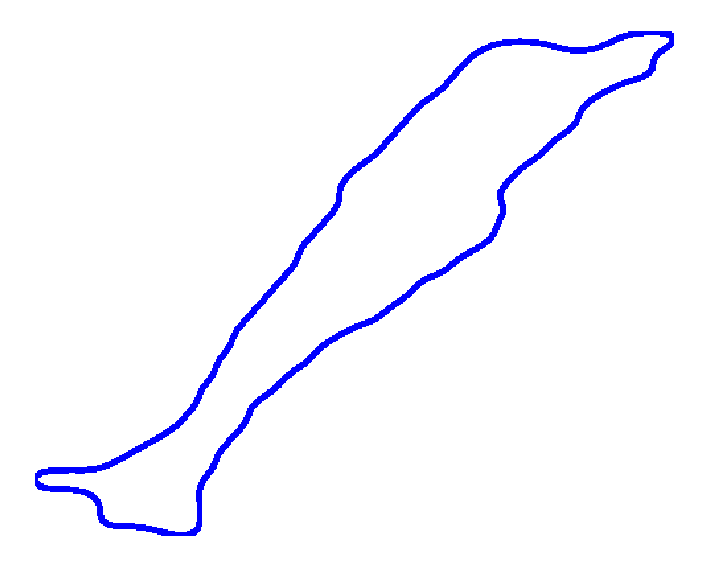}
&
\includegraphics[height=.075\textwidth]{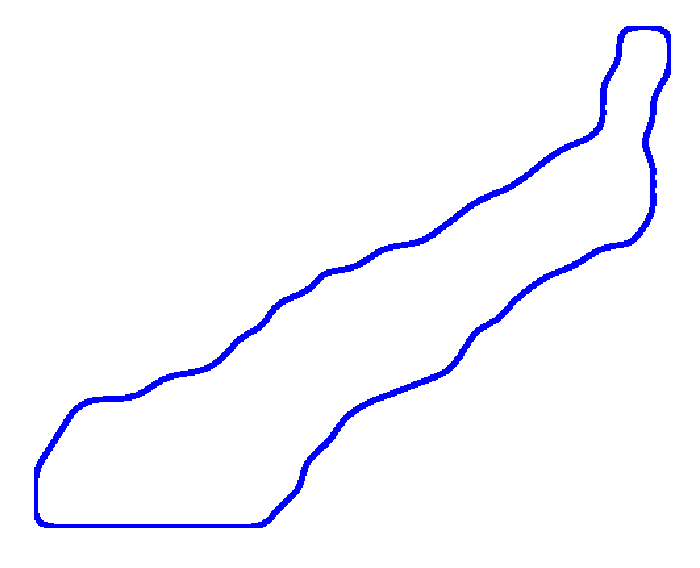}
&
\includegraphics[height=.075\textwidth]{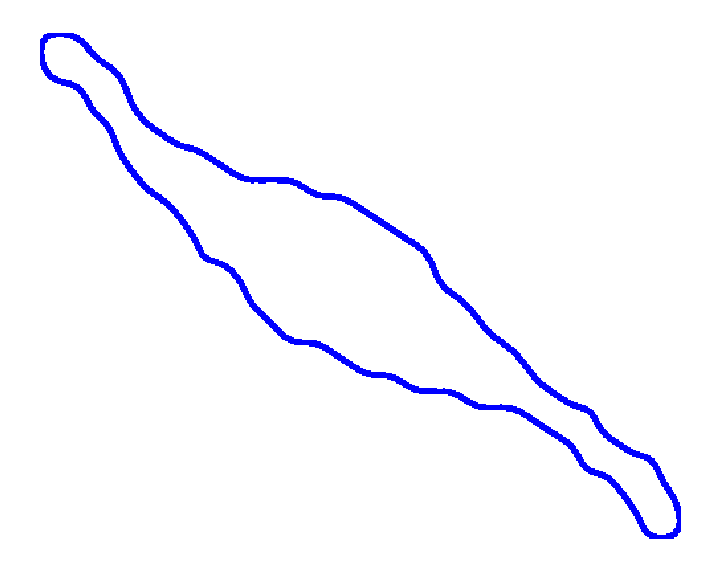}
&
\includegraphics[height=.075\textwidth]{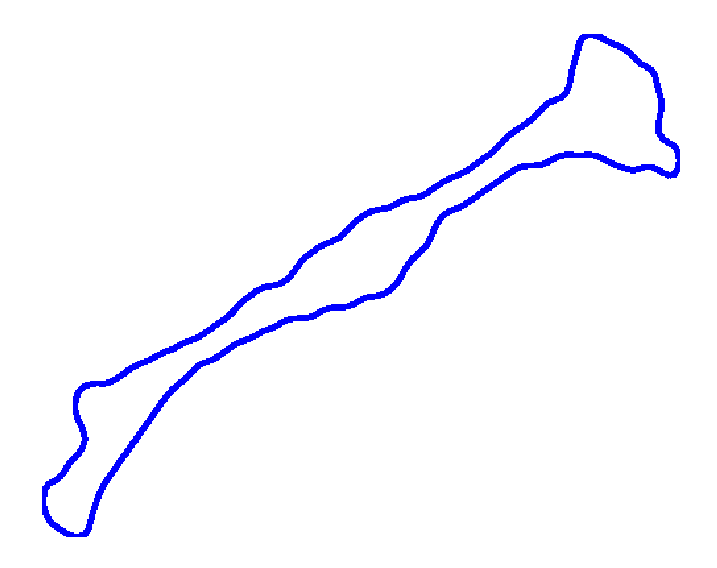}
&
\includegraphics[height=.075\textwidth]{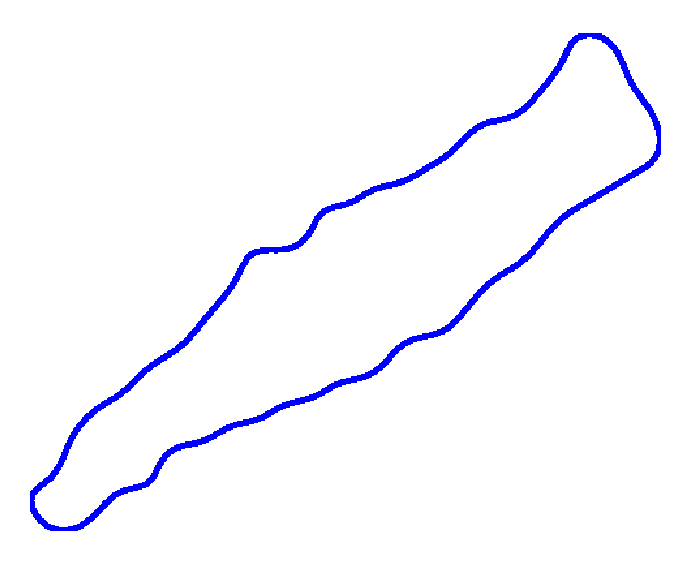}
\\ \hline
\end{tabular}
\end{center}
\caption{\label{F:example-curves}
Curve examples used in the experiments: synthetic curves (top row),
biological cell boundaries of type A and B (middle and bottom~rows).
}
\end{figure}

\begin{figure}
\begin{center}
\includegraphics[height=.4\textwidth]{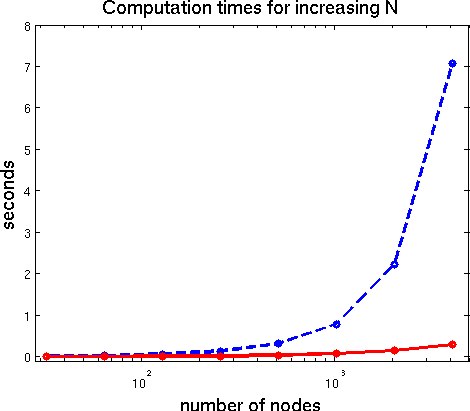}
\end{center}
\caption{\label{F:alignment-scalability}
Comparison of timings on an example curve,
between the slow $O(N^2)$ alignment algorithm (blue)
and the fast $O(N \log N)$ alignment algorithm proposed (red).
}
\end{figure}

We did an experimental validation of our algorithm by using different
closed curves as test cases. We used five synthetic curves (super-ellipse,
hippopede, bumps, lima\c{c}on, clover), five (biological) cell boundaries
of type A, and five (biological) cell boundaries of type B
(see Figure~\ref{F:example-curves}). For each of these reference curves,
we tried to align a template curve~to~it. We obtained the template curve by
changing the starting point of the reference curve from 0 to 0.25
and rotating the curve by $\frac{\pi}{3}$. We tested with the same uniform
arc length parametrization assigned to both curves.
We examined the running times and the alignment errors (quantified
by Equation~\eqref{E:disc-energy}) for increasing number of nodes $N$
on the curves. We observed consistent behaviour across different
cell examples. Our new $O(N \log N)$ algorithm was much faster
than the old $O(N^2)$ algorithm and had exactly the same alignment
error values that the old one had. These results are given
for~the~first~cell~of~type~B in Figure~\ref{F:alignment-scalability} and
Table~\ref{T:alignment-scalability}.
\begin{table}
\begin{center}
\begin{tabular}{|l|c|c|c|c|c|c|c|}
\hline
 & {\small N=64} & {\small 128} & {\small 256} & {\small 512} & {\small 1024} & {\small 2048} & {\small 4096} \\
\hline
Slow & .027 & .06 & .14 & .31 & .79 & 2.2 & 7.1 \\
Fast & .005 & .01 & .02 & .04 & .07 & .15 & .30 \\
Factor & 5.4x & 6.9x & 7.4x & 7.5x & 11x & 15x & 24x \\
\hline
Error & {\small 2e-2} & {\small 1e-2} & {\small 6e-3} &
       {\small 3e-3} & {\small 1.5e-3} & {\small 7e-4} & {\small 4e-4} \\
\hline
\end{tabular}
\end{center}
\smallskip
\caption{Timings (seconds) and errors for alignment algorithms.
The two algorithms produce same alignment errors.}
\label{T:alignment-scalability}
\end{table}

\section{Algorithm for elastic shape distances}\label{S:shape-dist}

In \cite{srivastava2}, Srivastava et al.\! proposed an elastic shape distance
between closed curves of unit length. In order to represent the shape of one
such curve $\beta$, $\beta: [0,1]\rightarrow \mathbb{R}^2$,
they defined a function called the
square root velocity (SRV) function of $\beta$ by
$q(t) = \dot{\beta}(t) / \|\dot{\beta}(t)\|^{1/2}$, \mbox{$t \in [0,1]$}, and an orbit
$[q]=\{\sqrt{\dot{\gamma}} R q(\gamma) \mid \gamma\in \Gamma, R\in{\it SO(2)}\}$
representing the shape of $\beta$ (encoding invariance with
respect to rotation and reparameterization by $R$ and $\gamma$, respectively;
translation invariance is taken care of automatically as
$\frac{d}{dt}(\beta(t)+c) = \dot{\beta}(t)$ for any constant planar point~$c$).
Here {\it SO(2)} is the set of all rotations in $\mathbb{R}^2$ and
$\Gamma$ is the set of all diffeomorphisms of $[0,1]$ into $[0,1]$, with
$\gamma(0)=0$, $\gamma(1)=1$, for $\gamma$ in~$\Gamma$.
We note that $\beta$ of unit length implies $\|q\|_{L^2}=1$ which in turn
implies $q$ is square-integrable. In this setting, the shape distance
between closed curves $\beta_1$ and $\beta_2$ of unit length is the
distance between the corresponding square-integrable optimally-matching SRV
functions $q^*_1(t) := q_1(t)$ and
$q^*_2(t) = \sqrt{\dot{\gamma}(t)} R q_2(t_0+\gamma(t))$, assumed to be
centered at the origin, where the triple $(t_0,R,\gamma)$ is the global
minimizer of mismatch energy
\begin{equation}\label{E:dist-energy}
E(t_0,R,\gamma) = \int_0^1 \| q_1(t)
- \sqrt{\dot{\gamma}(t)} R q_2(t_0 + \gamma(t)) \|^2 dt.
\end{equation}
Srivastava et al.\! minimize \eqref{E:dist-energy}
by looping over all $t_0$ candidates, computing optimal
$R$ with Kabsch algorithm for each $t_0$, and
then using dynamic programming to get best $\gamma$ for
each $(t_0,R)$ pair \cite{website,mio,srivastava2}.
This algorithm, which we call Approach~1,
is computationally expensive as it is $O(N^3)$, $N$ the number of nodes
per curve. Recently, an $O(N^2)$ iterative algorithm was proposed in
\cite{dogan} of observed subquadratic almost linear
time complexity, to optimize \eqref{E:dist-energy},
which we call Approach~2. It starts from an initial
triple of $(t_0,R,\gamma)$ ($(0,Id,\gamma(t)=t)$ in
\cite{dogan}),
and then updates $(t_0,R,\gamma)$ as it
alternates between optimizations with respect to $(t_0,R)$
and $\gamma$, until
\eqref{E:dist-energy} is minimized.
Here, assuming without loss of generality that $\beta_1$ and $\beta_2$ are
square-integrable and centered at the origin, we propose to use the
FFT-based optimal rigid alignment algorithm
of the previous section applied on $\beta_1(t)$ and $R \beta_2(t+t_0)$,
together with uniform resampling of $\beta_1$
and $\beta_2$ with respect to arc length, to compute optimal
$(t_0^*,R^*)$ of~\eqref{E:energy} and then start the algorithm in
\cite{dogan} with $(t_0^*,R^*,\gamma(t)=t)$.
Additionally, we propose to modify the algorithm in
\cite{dogan},
\emph{i.e.}, Approach~2,
as follows. The energy in \eqref{E:dist-energy} is
reformulated as
\begin{equation}\label{E:dist2-energy}
E(t_0,R,\gamma) = \int_0^1 \| R q_1(t+t_0)
- \sqrt{\dot{\gamma}(t)} q_2(\gamma(t)) \|^2 dt
\end{equation}
and each optimization of \eqref{E:dist2-energy} with respect to $(t_0,R)$
is then carried out using in the same way as above the FFT-based optimal rigid alignment
algorithm of the previous section (obviously without arc-length uniform resampling
of any of the functions involved)
applied now on square-integrable functions $q^*_1(t) := R q_1(t+t_0)$ and
$q^*_2(t) = \sqrt{\dot{\gamma}(t)} q_2(\gamma(t))$,
instead of $\beta_1(t)$ and $R \beta_2(t+t_0)$, with $q^*_1$~and $q^*_2$ assumed
to be centered at the origin.
Accordingly, as mentioned above, Approach~2 is initialized by pre-alignment,
\emph{i.e.}, starts from the initial triple $(t_0^*,R^*,\gamma(t)=t)$,
and then updates $(t_0,R,\gamma)$ as it
alternates between optimizations with respect to $(t_0,R)$ (using the FFT)
and $\gamma$ (using fast dynamic programming and nonlinear constrained
optimization as decribed in~\cite{dogan}), until
\eqref{E:dist2-energy} is minimized.
We found that modifying the algorithm in
\cite{dogan} with the FFT this way results in superb
performance in terms of computation times and computed minimizers
$(t_0,R,\gamma)$.
To evaluate our algorithm, we performed two sets of
experiments, first with synthetic curves, then with (biological) cell boundaries.

In the first set of experiments, we aimed to evaluate both
accuracy and scalability with respect to increasing~$N$.
For this, we took a synthetic curve $\beta_1$,
changed starting point from $0$ to $0.25$, rotated by $\frac{\pi}{3}$,
and reparameterized $\beta_2$ by $\gamma^1(t) = t+0.025\sin(4\pi t)$
and then $\gamma^2(t) = t + 1.6 t^2(t-1)^2$
(which induce mildly and strongly nonuniform distributions
of nodes, respectively).
The theoretical shape distance between $\beta_1$ and $\beta_2$ is zero.
We computed the shape distance using Approach 1, Approach 2
initialized by pre-alignment, but no FFT in the actual optimization,
then again Approach 2 now with the FFT used in the SRV function
optimization as well, for all of the five synthetic
curves in Figure~\ref{F:example-curves} with increasing~$N$.
The average shape distances and computation times are given
in Table~\ref{T:distances-scalability}. Approach 2 computed
zero for almost all distances, whereas Approach 1 computed significant
nonzero values. Moreover, the computation times of Approach 1
increased dramatically as we increased $N$, whereas the computation
time of Approach 2 grew very slowly, and was less than 1 s even
for curves with $N=1024$ nodes.

In the second set of experiments, we computed matrix of
pairwise distances for the 10 cell boundaries, each uniformly
sampled at $N=256$. Computing the $10\times 10$ distance
matrix took 9218 s with Approach 1, 125 s with Approach 2
and FFT-based optimization
(74X faster). The role of estimating $t_0$ and $R$
was not as apparent for this curve set. Still Approach 2
computed smaller distances than Approach 1 for 61 of the 100 pairs,
indicating better minimization of energy~\eqref{E:dist-energy}.
In fact, in only 6 of the remaining 39, Approach 1 computed
a significantly smaller distance than Approach 2
(difference $\ge 0.05$). If we look at a more relevant subset of the
cells, say type B curves, for which optimizing rotation plays a role,
we see that Approach 1 computed smaller distances in only 6 of the 25
pairs~(see~Table~\ref{T:cell-dist}).

\begin{table}
\begin{center}
\begin{tabular}{|l|c|c|c|c|c|c|}
\hline
{\it distance} & & {\small N=64} & {\small 128} & {\small 256} & {\small 512} & {\small 1024} \\
\hline
Approach 1 & $\gamma^{1}$ & .0184 & .0409 & .0362 & .0361 & N/A \\
           & $\gamma^{2}$ & .3344 & .3373 & .3445 & .3441 & N/A \\
\hline
Approach 2 & $\gamma^{1,2}$ & .0037 & 0 & 0 & 0 & 0 \\
\hline
\end{tabular}
\begin{tabular}{|l|c|c|c|c|c|c|}
\hline
{\it timing} & & {\small N=64} & {\small 128} & {\small 256} & {\small 512} & {\small 1024} \\
\hline
Approach 1 & $\gamma^{1,2}$ & 1.45 & 11.5  & 94 & 756 & N/A \\
\hline
Approach 2 & $\gamma^{1}$ & 0.20 & 0.20 & 0.33 & 0.63 & 1.75 \\
\ (no FFT) & $\gamma^{2}$ & 0.34 & 0.29 & 0.37 & 0.70 & 1.68 \\
\hline
Approach 2 & $\gamma^{1}$ & 0.15 & 0.14 & 0.19 & 0.35 & 0.85 \\
\ (FFT)    & $\gamma^{2}$ & 0.16 & 0.19 & 0.26 & 0.42 & 0.96 \\
\hline
\end{tabular}
\end{center}
\caption{Average shape distances and timings (in seconds)
for the synthetic curves shown in Figure~\ref{F:example-curves},
using Approach 1, Approach 2 without FFT speed-up,
Approach 2 with it.}
\label{T:distances-scalability}
\end{table}

\begin{table}
\begin{center}
\begin{tabular}{|c|c|c|c|c|}
\hline
 0.0/0.0  & .406/.373 & .363/.343 & .316/.303 & {\bf .351/.378} \\ \hline
{\bf .386/.409} &  0.0/0.0  & .430/.414 & .358/.337 & .376/.365 \\ \hline
{\bf .343/.372} & .434/.421 &  0.0/0.0  & .375/.344 & .317/.295 \\ \hline
{\bf .314/.333} & .367/.363 & {\bf .384/.395} &  0.0/0.0  & .327/.316 \\ \hline
.349/.345 & .376/.372 & {\bf .306/.309} & .328/.305 &  0.0/0.0 \\
\hline
\end{tabular}
\end{center}
\caption{Matrix of pairwise shape distances of type B cells.
The first and second values of a pair computed
by Approach 1 and Approach 2 with FFT speed-up, respectively.
The instances when Approach 1 computed smaller distances than
Approach 2 shown in bold.}
\label{T:cell-dist}
\end{table}

\section{Conclusions}\label{S:conclusion}

In this paper, we propose a fast curve rigid alignment algorithm of
$O(N \log N)$ time complexity based on the FFT, $N$~the number of nodes
per curve. Given two closed curves, together with uniform
resampling of the curves with respect to arc length, our algorithm
computes the optimal starting point $t_0$ and rotation $R$ that must
be applied on one curve in order to rigidly align it with the other curve.
This is the main contribution of this paper.
Additionally, we describe how this fast rigid alignment algorithm can be
used to initialize the iterative algorithm in~\cite{dogan}
for computing elastic shape distances, and how it can be
used for the $(t_0,R)$ optimization step of the same algorithm.
The resulting new algorithm computes then accurate shape distances at
a fraction of the cost of the original one, and scales
very well to curves with large numbers of nodes.


\begin{thebibliography}{5}

\bibitem{website}
Code from {S}tatistical {S}hape {A}nalysis and {M}odeling {G}roup, {F}lorida
  {S}tate {U}niversity.
http://ssamg.stat.fsu.edu/downloads/ClosedCurves2D3D.zip.
Accessed: 2014-06-20.

\bibitem{ayache}
N.~J. Ayache and O.~D. Faugeras.
\newblock Hyper: A new approach for the recognition and positioning of
  two-dimensional objects.
\newblock \emph{IEEE Transactions on Pattern Analysis and Machine
  Intelligence}, 8\penalty0 (1):\penalty0 44--54, 1986.

\bibitem{cohen}
I.~Cohen, N.~J. Ayache, and P.~Sulger.
\newblock Tracking points on deformable objects using curvature information.
\newblock In \emph{Proceedings European Conference on Computer Vision}, pages
  136--144, 1993.

\bibitem{cui}
M.~Cui, J.~Femiani, J.~Hu, P.~Wonka, and A.~Razdan.
\newblock Curve matching for open 2d curves.
\newblock \emph{Pattern Recognition Letters}, 30\penalty0 (1):\penalty0 1--10,
  2009.

\bibitem{dogan}
G.~Do\u{g}an, J.~Bernal, and C.R. Hagwood.
\newblock A fast algorithm for elastic shape distances between closed planar
  curves.
\newblock In \emph{Proceedings of the IEEE Conference on Computer Vision and
  Pattern Recognition}, June 2015.

\bibitem{kabsch1}
W.~Kabsch.
\newblock A solution for the best rotation to relate two sets of vectors.
\newblock \emph{Acta Crystallographica Section A: Crystal Physics}, 32\penalty0
  (5):\penalty0 922--923, 1976.

\bibitem{kabsch2}
W.~Kabsch.
\newblock A discussion of the solution for the best rotation to relate two sets
  of vectors.
\newblock \emph{Acta Crystallographica Section A: Crystal Physics}, 34\penalty0
  (5):\penalty0 827--828, 1978.

\bibitem{larsen}
R.~Larsen.
\newblock L1 generalized procrustes 2d shape alignment.
\newblock \emph{Journal of Mathematical Imaging and Vision}, 31\penalty0
  (2-3):\penalty0 189--194, 2008.

\bibitem{li}
H.~Li, T.~Shen, and X.~Huang.
\newblock Approximately global optimization for robust alignment of generalized
  shapes.
\newblock \emph{Pattern Analysis and Machine Intelligence, IEEE Transactions
  on}, 33\penalty0 (6):\penalty0 1116--1131, June 2011.

\bibitem{mio}
W.~Mio, A.~Srivastava, and S.~Joshi.
\newblock On shape of plane elastic curves.
\newblock \emph{International Journal of Computer Vision}, 73\penalty0
  (3):\penalty0 307--324, 2007.

\bibitem{schwartz}
J.~T. Schwartz and M.~Sharir.
\newblock Identification of partially obscured objects in two and three
  dimensions by matching noisy characteristic curves.
\newblock \emph{Int'l J. Robotics Research}, 6\penalty0 (2):\penalty0 29--44,
  1987.

\bibitem{sebastian}
T.~B. Sebastian, P.~N. Klein, and B.~B. Kimia.
\newblock On aligning curves.
\newblock \emph{IEEE Transactions on Pattern Analysis and Machine
  Intelligence}, 25\penalty0 (1):\penalty0 116--125, 2003.

\bibitem{srivastava2}
A.~Srivastava, E.~Klassen, S.H. Joshi, and I.H. Jermyn.
\newblock Shape analysis of elastic curves in {E}uclidean spaces.
\newblock \emph{Pattern Analysis and Machine Intelligence, IEEE Transactions
  on}, 33\penalty0 (7):\penalty0 1415--1428, 2011.

\bibitem{umeyama}
S.~Umeyama.
\newblock Parameterized point pattern matching and its application to
  recognition of object families.
\newblock \emph{IEEE Transactions on Pattern Analysis and Machine
  Intelligence}, 15\penalty0 (2):\penalty0 136--144, 1993.

\end{thebibliography}
\end{document}